\documentclass{amsart}
\usepackage{graphics}

\newcommand{\II}{\mathbb I}
\newcommand{\ID}{\mathbb D}
\newcommand{\IT}{\mathbb T}
\newcommand{\IC}{\mathbb C}
\newcommand{\IZ}{\mathbb Z}

\newcommand{\IR}{\mathbb R}
\newcommand{\IN}{\mathbb N}
\newcommand{\M}{\mathcal M}
\newcommand{\U}{\mathcal U}
\newcommand{\IH}{\mathbb H}
\newcommand{\id}{\mathrm{id}}
\newcommand{\e}{\varepsilon}

\newcommand{\non}{\mathrm{non}}
\newcommand{\mesh}{\mathrm{mesh}}
\newcommand{\asdim}{\mathrm{asdim}}

\newcommand{\ms}{\mathsf{ms}}
\newcommand{\bms}{\mathsf{bms}}
\newcommand{\dom}{\mathrm{dom}}
\newcommand{\tms}{\mathsf{m}\tilde{\mathsf s}}

\newcommand{\Iso}{\mathrm{Iso}}
\newcommand{\MS}{\mathsf{MS}}

\newcommand{\tS}{\widetilde{\mathcal S}}
\newcommand{\mds}{\mathsf{mds}}

\newcommand{\card}{\mathrm{card}}
\newcommand{\Borel}{\mathsf{Borel}}
\newcommand{\diam}{\mathrm{diam}}
\newcommand{\cp}{\mathsf{c}}
\newcommand{\Bcp}{\mathsf{Bc}}
\newcommand{\supp}{\mathrm{supp}}

\newtheorem{theorem}{Theorem}[section]
\newtheorem{problem}[theorem]{Problem}
\newtheorem{question}[theorem]{Question}
\newtheorem{proposition}[theorem]{Proposition}
\newtheorem{lemma}[theorem]{Lemma}
\newtheorem{corollary}[theorem]{Theorem}
\newtheorem{example}[theorem]{Example}

\theoremstyle{definition}
\newtheorem{definition}{Definition}

\title[Symmetry and Colorings]{Symmetry and colorings:\\ some results and open problems, II}
\author{Taras Banakh}
\address{Ivan Franko University of Lviv (Ukraine) and Jan Kochanowski University in Kielce (Poland)}
\email{tbanakh@yahoo.com}

\begin{document}
\begin{abstract} We survey some principal results and open problems related to colorings of geometric and algebraic objects endowed with symmetries. We concentrate the exposition on the maximal symmetry numbers $\ms_r(\mathbf X)$ of such objects.
\end{abstract}
\maketitle

\section{Introduction}

In this paper we give a survey of principal results and open problems related to colorings and symmetry. The paper can be considered as an update of the earlier surveys \cite{BP01} and \cite{BVV00e}, written more than 10 years ago. It can also serve as lecture notes of the cycle of 16 lectures delivered by the author in Silesian University (Poland) in the fall Semester of 2011. The paper contains no proofs. Our aim is to describe the general picture of the theory and refer to original papers for proofs of mentioned results.

The subject discussed in this paper traces its history back to S.Sidon who posed in 1931 the
following problem related to the Additive Number Theory \cite{PS95}: given a  subset $A\subset\IN$ and a number $n\in\IN$ estimate the quantity of
solutions $(x,y)\in A\times A$ of the equation $x+y=n$. Observe that this is essentially a question about symmetric subsets of $A$: it is easy to see that the set $S(n)=\{x\in A:n-x\in A\}$ is symmetric with respect to the point $n/2$ and its cardinality is just equal to the quantity of pairs $(x,y)\in A\times A$ with $x+y=n$. Thus the Sidon question can be
reformulated in terms of symmetry: given a subset $A\subset \IN$ estimate the size of a maximal symmetric subset of $A$. This observation is of crucial importance since symmetric subsets can be naturally defined for many algebraic and geometric objects.

We shall treat symmetry using the methods of Ramsey Theory. Ramsey Theory is a science about pieces of order in chaos. A typical Ramsey result says that for any partition of certain highly organized object $X$ into pieces one of the pieces still contains a relatively large subobject with high organization, see \cite{GRS}. Instead of partitions into $r$ pieces, it is more convenient to speak about colorings of $X$ in $r$ colors. Then subsets that lie in the pieces of the partition are colored by the same color and are called {\em monochromatic}.

A typical example of a Ramsey-theoretic result is the Van der Waerden Theorem saying that for any finite partition $\IN=A_1\cup\dots\cup A_r$ of the set $\IN$, one of the pieces of the partition contains arbitrarily long arithmetic progressions, see \cite{GRS}. In the chromatic terminology this theorem says that for any finite coloring of $\IN$ there are arbitrarily long monochromatic arithmetic progressions.

In this paper we shall treat the symmetry by Ramsey-theoretic methods and for colorings of geometric or algebraic objects $X$, shall look for ``large'' monochromatic subsets which are symmetric (so, in a sense have some inner structure). The ``largeness'' of subsets of such object $X$ will be evaluated with help of a ``measure'' $\mu:\dom(\mu)\to[\mu(\emptyset),\mu(X)]$ defined on a $\sigma$-algebra $\dom(\mu)$ of subsets of $X$, while symmetric subsets will be defined by fixing a family $\mathcal S$ of measure-preserving bijective transformations of $X$, called admissible symmetries of $X$.
For such a triple $\mathbf X=(X,\mu,\mathcal S)$ (called a symmetry measure space) we shall evaluate the maximal symmetry numbers $\ms_r(\mathbf X)$ equal to the maximal measure of a symmetric subset that can be found in each measurable $r$-coloring of $X$. In subsequent sections we shall calculate or evaluate this maximal symmetry number $\ms_r(\mathbf X)$ for various symmetry measure spaces $\mathbf X=(X,\mu,\mathcal S)$ that naturally appear in Geometry and Algebra.

\section{Symmetry Measure spaces}

In this section we introduce the notion of a symmetry measure space $\mathbf X=(X,\mu,\mathcal S)$.
This is a triple consisting of a set $X$, a measure $\mu$ on $X$ and a family $\mathcal S$ of bijective measure-preserving transformations of $X$, called  admissible symmetries of $X$.

By a {\em measure} on $X$ we understand any monotone function $\mu:\dom(\mu)\to[\mu(\emptyset),\mu(X)]$ defined on a $\sigma$-algebra $\dom(\mu)$ of subsets of $X$ and taking values in a complete linearly ordered set $[\mu(\emptyset),\mu(X)]$ with minimal element $\mu(\emptyset)$ and maximal element $\mu(X)$. The completeness of $[\mu(\emptyset),\mu(X)]$ means that each subset $A\subset [\mu(\emptyset),\mu(X)]$ has the exact lower and upper bounds $\inf(A)$ and $\sup(A)$.
The monotonicity of $\mu$ means that $\mu(A)\le\mu(B)$ for any subsets $A\subset B$ that belong to the $\sigma$-algebra $\dom(\mu)$. The pair $(X,\mu)$ will be called a {\em measure space}. Elements of the $\sigma$-algebra $\dom(\mu)$ are called {\em $\mu$-measurable} (or just {\em measurable}) subsets of the measure space $(X,\mu)$.

The other component of a symmetry measure space is a family $\mathcal S$ of bijective transformations of $X$, called {\em admissible symmetries} (or just {\em symmetries}) of $X$.
We assume that the family $\mathcal S$ is symmetric in the sense that $s^{-1}\in\mathcal S$ for each $\mathcal S$, and that each admissible symmetry $s\in S$ is {\em measure-preserving} in the sense that for any subset $A\in\dom(\mu)$ its image $s(A)$ belongs to the $\sigma$-algebra $\dom(\mu)$ and has measure $\mu(s(A))=\mu(A)$. A symmetry $s:X\to X$ is called {\em involutive} if $s^{-1}=s$, i.e., $s\circ s$ coincides with the identity transformation $\id_X$ of $X$.

A subset $A\subset X$ of a symmetry measure space $(X,\mu,\mathcal S)$ is called {\em symmetric} (more precisely, {\em $\mathcal S$-symmetric}) if $A=s(A)$ for some symmetry $s\in\mathcal S$. This notion is trivial if $\mathcal S$ contains the identity transformation of $X$ (in which case, each subset of $X$ is $\mathcal S$-symmetric). So, as a rule we shall assume that $\mathcal S$ does not contain the identity transformation $\id_X$ of $X$, but the case $\mathcal S=\{\id_X\}$ is not completely excluded.

Observe that for any set $A\subset X$ and an admissible symmetry $s$ the intersection
$$\bigcap_{n\in\IZ}s^n(A)$$ is the largest $\{s\}$-symmetric subset of $A$.
If the symmetry $s$ is involutive, then this intersection coincides with $A\cap s(A)$. If the set $A$ belongs to the $\sigma$-algebra $\dom(\mu)$ of $\mu$-measurable sets then the set $\bigcap_{n\in\IZ}s^n(A)$ also belongs to $\dom(\mu)$.

\section{The maximal symmetry numbers $\ms_r(\mathbf X)$ of a symmetry measure space $\mathbf X$}

In this section, for each symmetry measure space $\mathbf X=(X,\mu,\mathcal S)$ we define its maximal symmetry numbers $\ms_r(\mathbf X)$ parametrized by cardinal numbers $r$.

First we define the maximal symmetry number $\ms(A)$ of a measurable subset $A$ of $X$.
By definition,
$$\ms(A)=\sup\{\mu(B):\mbox{$B$ is a $\mu$-measurable $\mathcal S$-symmetric  subset of $A$}\}$$is the maximal measure of a $\mu$-measurable $\mathcal S$-symmetric subset of $A$.
It is clear that
$$\ms(A)=\sup_{s\in\mathcal S}\mu\big(\textstyle{\bigcap\limits_{n\in\IZ}}s^n(A)\big),$$
where $\bigcap\limits_{n\in\IZ}s^n(A)$ is the maximal $\{s\}$-symmetric subset of $A$.
If each admissible symmetry $s\in \mathcal S$ is involutive, then
$$\ms(A)=\sup_{s\in\mathcal S}\mu(A\cap s(A)).$$


Next, for every cardinal number $r$ put
$$\ms_r(A)=\inf_{\chi:A\to r}\sup_{i\in r}\ms(\chi^{-1}(i)),$$
where the infimum is taken over all $\mu$-measurable colorings $\chi:A\to r$ of $A$. The $\mu$-measurability of $\chi$ means that each color class $\chi^{-1}(i)$, $i\in r$, belongs to the $\sigma$-algebra $\dom(\mu)$ of $\mu$-measurable sets. The maximal symmetry number $\ms_r(A)$ takes its values in the linearly ordered set $[\mu(\emptyset),\mu(X)]$.

If we want to specify the measure $\mu$ and the family of symmetries $\mathcal S$, then  we shall use the expanded form $\ms_r(A,\mu,\mathcal S)$ instead of $\ms_r(A)$. Also for a symmetry measure space $\mathbf X=(X,\mu,\mathcal S)$ the maximal symmetry numbers $\ms_r(X,\mu,\mathcal S)$ will be denoted by $\ms_r(\mathbf X)$.

The definition of the number $\ms_r(A)$ implies that
\begin{itemize}
\item for any $a<\ms_r(A)$ and each measurable $r$-coloring of $A$ there is a monochrome $\mathcal S$-symmetric measurable subset $B\subset A$ with measure $\mu(B)>a$;
\item for any $b>\ms_r(A)$ there is a measurable $r$-coloring of $A$ without monochrome $\mathcal S$-symmetric measurable subset $B\subset A$ of measure $\mu(B)\ge b$.
\end{itemize}

The maximal symmetry number $\ms_r(A,\mu,\mathcal S)$ is monotone with respect to $r$, $A$, and $\mathcal S$ in the sense that
$$\ms_{r}(A,\mu,\mathcal S)\le\ms_{r'}(A',\mu,\mathcal S')$$ for any
cardinals $r\ge r'$, measurable sets $A\subset A'$, and families $\mathcal S\subset\mathcal S'$ of admissible symmetries.

This inequality can be also derived from the following general fact:

\begin{proposition} For two symmetry measure spaces $\mathbf X=(X,\mu_X,\mathcal S_X)$ and $\mathbf Y=(Y,\mu_Y,\mathcal S_Y)$ the inequality $\ms_r(\mathbf X)\le\ms_r(\mathbf Y)$ holds for every cardinal $r$ provided there is a function $f:X\to Y$ such that
\begin{itemize}
\item for any $\mu_Y$-measurable subset $B\subset Y$ the preimage $f^{-1}(B)$ is $\mu_X$-measurable and has measure $\mu_X(f^{-1}(B))\le \mu_Y(B)$;
\item for any admissible symmetry $s_X\in\mathcal S_X$ there is an admissible symmetry $s_Y\in\mathcal S_Y$ such that $s_Y\circ f=f\circ s_X$.
\end{itemize}
\end{proposition}

\section{The maximal symmetry number of some geometric objects}

In this section we shall evaluate the maximal symmetry numbers $\ms_r(\mathbf X)$ of some symmetry measure spaces $\mathbf X$ that naturally appear in Geometry. Among typical spaces of this sort let us mention spheres and balls in Euclidean spaces. We shall give upper and lower bound for such symmetry measure spaces having geometric nature.

The following general theorem proved in \cite{BVV00i} gives a typical upper bound for the maximal symmetry numbers $\ms_r(\mathbf X)$ of compact subsets of connected Riemannian manifolds, endowed with their canonical measure and canonical metric.

\begin{theorem}[{\cite{BVV00i}}]\label{t4.1} Let $M$ be a Riemannian manifold, endowed with its canonical Borel measure $\mu$ and the canonical metric $d$. Let $X$ be a compact subset of $M$ and $\mathcal S$ be the family of all non-identity isometries $f:M\to M$ such that $f(X)=X$. Then
$$\ms_r(X,\mu,\mathcal S)\le \frac{\mu(X)}{r^2}\mbox{ \ for all $r\in\IN$},$$ which means that for every $\e>0$ there is a measurable $r$-coloring of $X$ without monochromatic $\mathcal S$-symmetric subsets of measure $\ge \frac{\mu(X)}{r^2}+\e$.
 \end{theorem}

The lower bound $\frac{\mu(X)}{r^2}\le\ms_r(\mathbf X)$ can be proved for symmetry measure spaces $\mathbf X=(X,\mu,\mathcal S)$ with sufficiently rich family of admissible symmetries, in particular, for spaces having rotational or spherical symmetries.

To define spaces with such symmetries, we first introduce the symmetry measure structure on $n$-dimensional spheres $S^n$ and cylinders $S^n\times[0,a]$. We endow the $n$-dimensional sphere
$S^n=\{x\in\IR^{n+1}:\|x\|=1\}$ with a unique probability Borel measure $\mu$, invariant under the action of the group $O(n+1)$ of orthogonal transformations of the Euclidean space $\IR^{n+1}$.
For the family of admissible symmetries there are many natural choices. The largest is the family $\mathcal I\setminus\{\id\}$ of all non-identity isometries of $S^n$. This family contains the subfamily $\mathcal R$ consisting of all reflections of $S^n$ with respect to the hyperplanes passing through the center of the sphere.

\begin{theorem}[{\cite{BVV00i}}]\label{t4.2} For every $r,n\in\IN$ we get $$\ms_r(S^n,\mathcal R)=\ms_r(S^n,\mathcal I\setminus\{\id\})=\frac1{r^2}.$$ Moreover, for each measurable coloring of $S^n$ by $r\ge 2$ colors there is a monochromatic $\mathcal R$-symmetric subset $A\subset S^n$ of measure $\mu(A)>\frac1{r^2}$.
\end{theorem}

A similar result is true for cylinders $S^n\times(a,b)$ over spheres, and suitable images of such cylinders. Here we endow the interval $(a,b)$ with the Lebesgue measure and the trivial family $\{\id\}$ of admissible symmetries.

The cylinder $S^n\times(a,b)$ will be considered as a symmetry measure space endowed with the product measure $\lambda$ and the family $\mathcal R\times\{\id\}$ of admissible symmetries. So, each admissible symmetry of $S^n\times(a,b)$ is of the form $$s\times\id:(x,t)\mapsto (s(x),t)$$for some reflection $s\in\mathcal R$. Geometrically, this means that admissible symmetries of the cylinder $S^n\times(a,b)$ are reflections with respect to the hyperplanes in $\IR^{n+1}\times\IR$ that contain the axis $\{0\}\times\IR$.

We say that a symmetry measure space $\mathbf X=(X,\mu,\mathcal S)$ has {\em $S^n$-symmetry} if the measure $\mu$ takes its values in the real line $\IR$ and there is a function $f:S^n\times(0,\mu(X))\to X$ such that
\begin{enumerate}
\item for any $\mu$-measurable subset $B\subset X$ the set $f^{-1}(B)$ is $\lambda$-measurable and has measure $\lambda(f^{-1}(B))=\mu(B)$;
\item for any reflection $r\in\mathcal R\times\{\id\}$ of the cylinder $S^n\times(0,\mu(X))$ there is an admissible symmetry $s\in\mathcal S$ such that $s\circ f=f\circ r$.
\end{enumerate}
These conditions imply that $\mu(\emptyset)=\lambda(f^{-1}(\emptyset))=0$, so $[\mu(\emptyset),\mu(X)]=[0,\mu(X)]\subset\IR$.

We shall say that a symmetry measure space $\mathbf X$ has {\em rotation symmetry} (resp. {\em spherical symmetry}) if $\mathbf X$ has $S^1$-symmetry (resp. $S^2$-symmetry).

Since the sphere $S^2$ has rotational symmetry, each symmetry measure space with spherical symmetry has rotation symmetry. Observe that all spheres and balls in the spaces $\IR^n$, $n\ge 3$, have spherical symmetry.

\begin{theorem}[{\cite{BVV00i}}] If a symmetry measure space $\mathbf X=(X,\mu,\mathcal S)$ has rotation symmetry, then $\ms_r(\mathbf X)\ge \frac{\mu(X)}{r^2}$ for every $r\in\IN$.
\end{theorem}

This theorem can be derived from the following more precise fact, proved in \cite{BVV00i} and \cite{BVV00e}:

\begin{proposition}\label{p4.4} If a symmetry measure space $\mathbf X=(X,\mu,\mathcal S)$ with $\mu(X)=1$ has:
\begin{itemize}
\item rotation symmetry, then each $\mu$-measurable subset $A\subset X$ contains an $\mathcal S$-symmetric subset $B\subset A$ of measure $\mu(B)\ge \mu(A)^2$;
\item spherical symmetry, then each $\mu$-measurable subset $A\subset X$ of measure $0<\mu(A)<1$ contains an $\mathcal S$-symmetric subset $B\subset A$ of measure  $\mu(B)>\mu(A)^2$.
\end{itemize}
\end{proposition}

The circle $$S^1=\{z\in\IC:|z|=1\}=\{e^{i\varphi}:0\le\varphi<2\pi\}$$ on the complex plane has no spherical symmetry but by Theorem~\ref{t4.2} still has the stronger  Ramsey property stated in Proposition~\ref{p4.4}: each subset $A\subset S^1$ of measure $0<A<1=\mu(S^1)$ contains an $\mathcal R$-symmetric subset $B\subset A$ of measure $\mu(B)>\mu(A)^2$.

In contrast, the 2-dimensional ball $B^2\subset\IR^2$ and the cylinder $S^1\times[0,1]$ do not have this stronger property. A counterexample is given by the subset
$$A=\big\{(e^{i\varphi},t)\in S^1\times [0,1]:\pi t\le\varphi<\pi+\pi t\big\}$$
of measure $\frac12$ which contains no symmetric subset of measure $>\frac14$.
Applying to $A$ the measure-preserving transformation $$f:S^1\times [0,1]\to \ID,\;\;\;f:(e^{i\varphi},t)\mapsto e^{i\varphi}\sqrt{t/\pi}$$of the cylinder $S^1\times[0,1]$ onto the disk $\ID=\{z\in\IC:|z|\le1/\sqrt{\pi}\}$ of unit area, we get the subset
$$F=\{r e^{i\varphi}:r^2\pi^2\le\varphi<r^2\pi^2+\pi\}$$ of the disc $\ID$ bounded by the Fermat spiral $\varphi=r^2\pi^2$. The set $F$ has area $\frac12$ but contains no axially symmetric subset of measure $>\frac14$. By its form it resembles one of the congruent pieces in the famous Yin-Yang symbol:

\centerline{\includegraphics{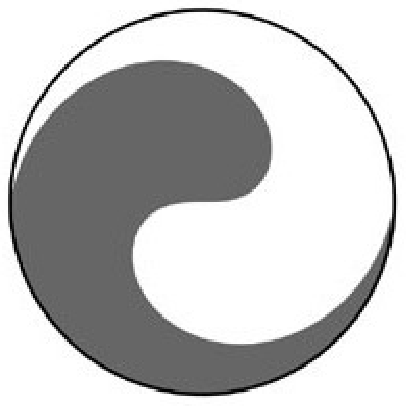}}

This property of $F$ was used in the paper \cite{BVV10} devoted to finding the canonical form of the Yin-Yang symbol.

A measurable subset $A$ of the disk $\ID=\{z\in\IC:\pi|z|^2\le 1\}$ of unit area will be called {\em equilibrial} if $A$ has area $\lambda(A)=\frac12$ and for any axial symmetry $s:\ID\to \ID$ of the disc the maximal $\{s\}$-symmetric subset $A\cap s(A)$ of $A$ has area $\lambda(A\cap s(A))=\frac14$.

\begin{definition} A curve $\beta$ in the disc $\ID$ will be called a {\em Yin-Yang curve} if it has the following 5 properties:
\begin{enumerate}
\item $\beta$ separates $\ID$ into two congruent equilibrial pieces $A$ and $-A$;
\item $\beta$ meets each concentric circle $C_a=\{z\in\IC:|z|=a\}$, $0<a\le1/\sqrt{\pi}$, in two opposite points;
\item $\beta$ has one-point intersection with each radius $R_\varphi=\{re^{i\varphi}:0<r\le1/\sqrt{\pi}\}$, $\varphi\in[0,2\pi)$;
\item the curve $\beta$ is smooth (i.e., has a regular $C^\infty$-parametrization);
\item the curve $\beta$ is algebraic in polar coordinates (i.e., its points satisfy a polynomial equation $P(r,\varphi)=0$).
\end{enumerate}
\end{definition}

\begin{theorem}[{\cite{BVV10}}] Up to the congruence the Fermat spiral $\beta=\{re^{i\varphi}\in\ID:\pi^2r^2=\varphi\}$ is a unique Yin-Yang curve.
\end{theorem}

\section{Maximal symmetry numbers $\ms_r(X)$ of $G$-spaces}

In this section we study the maximal symmetry numbers $\ms_r(X)$ for a measure space $(X,\mu)$ endowed with a measure-preserving action of a group $G$.

By an action of $G$ on $X$ we understand a function $\alpha:G\times X\to X$, $\alpha:(g,x)\mapsto gx$, such that
\begin{itemize}
\item $g(hx)=(gh)x$ for any $g,h\in G$ and $x\in X$;
\item $1x=x$ for all $x\in X$.
\end{itemize}
Here by $1$ we denote the neutral element of the group $G$.

An action $\alpha:G\times X\to X$ of a group $G$ on a measure space $(X,\mu)$ is called {\em measure-preserving} if for any element $g\in G$ the left shift $l_g:X\to X$, $l_g:x\mapsto gx$, is measure-preserving in the sense that for any $\mu$-measurable subset $A\subset X$ its shift $gA$ is $\mu$-measurable and has measure $\mu(gA)=\mu(A)$.

A measure space $X$ endowed with a measure-preserving action $\alpha$ of a group $G$ will be called a {\em measure $G$-space}.

Each measure $G$-space $(X,\mu,\alpha)$ will be considered as a symmetry measure space $\mathbf X=(X,\mu,\mathcal S)$ whose family of admissible symmetries $\mathcal S=\big\{l_g:g\in G\setminus\{1\}\big\}$ consists of  left shifts by non-identity elements of the group $G$. By $\ms_r(X)$ we shall denote the maximal symmetry number of that symmetry measure space.

Below $\II$ denotes the interval $[0,1]$ endowed with the standard Lebesgue
measure $\lambda$. For every compact topological group $G$ we consider the cylinder
$G\times \II$ as a left $G$-space (with the action $g\cdot(h,t)=(gh,t)$ for $g\in
G$ and $(h,t)\in G\times \II$) endowed with the product measure $\mu\times
\lambda$, where $\mu$ is the Haar measure on $G$.

\begin{theorem}[{\cite{BP01}}]\label{t5.1} Let $X$ be a compact metrizable space endowed with a $\sigma$-additive Borel probability measure $\mu$ and a continuous measure-preserving action of a finite group $G$. The for every finite number $r$,
\begin{enumerate}
\item $\ms_r(X)\ge \ms_r(G\times \II)$;
\item $\ms_r(X)=\ms_r(G\times \II,r)$ if the measure $\mu$ is atomless and
$|G\cdot x|=|G|$ for $\mu$-almost all $x\in X$.
\end{enumerate}
\end{theorem}

Thus the problem of calculating $\ms_r(X)$ reduces to calculating
the maximal symmetry numbers $\ms_r(G\times \II)$ of cylinder over finite groups $G$. Below $\mathcal D_{2n}$ is
a dihedral group, i.e., the group of all isometries of the regular
$n$-gon. For $n=2$ the dihedral group $\mathcal D_{2n}$ coincides with the Kleinian group
(isomorphic to $\IZ_2\times\IZ_2$).

The following unexpected result belongs to Ya.Vorobets and is taken from \cite{BP01}.

\begin{theorem}\label{t5.2}
\begin{enumerate}
\item For any subgroup $H$ of a finite group $G$ we have $\ms_r(H\times \II)\le
\ms_r(G\times \II)$ for every $r\in\IN$.
\item $\ms_r(\IZ_n\times \II)=0$ for every $n,r>1$.
\item For every prime number $p$
\begin{itemize}
\item  $\ms_2(\mathcal D_{2p}\times \II)=\frac{p-1}{p^2+2p-2}$,
\item $\ms_3(\mathcal D_{2p}\times \II)=\frac1{3p^2+6}$ and
\item $\ms_r(\mathcal D_{2p})=0$ if $r\ge 4$.
\end{itemize}
\item For any $n,r>1$ \ $\ms_r(\mathcal D_{2n})=\ms_r(\mathcal D_{2p})$, where $p$
is the minimal prime number dividing $n$;
\item If $r=2^k$ for some $k\ge 0$, then
$\ms_r(\IZ_2^n\times I)=\frac1{r^2}\cdot\frac{2^n-r}{2^n-1}$ for every $n>k$.
\end{enumerate}
\end{theorem}

Theorems~\ref{t5.1} and \ref{t5.2} allow us to classify the maximal symmetry numbers
$\ms_r(X)$ for all compact
subsets $X$ of positive Lebesgue measure in the plane, endowed with the family $\mathcal S=\Iso(X)\setminus\{\id\}$ of all nonidentity isometries of $X$.

\begin{theorem}[{\cite{BP01}}]\label{t5.3} For every compact subset $X\subset\IR^2$ of positive
Lebesgue measure  and  every $r\ge 2$ we have
$$
\ms_r(X)=\begin{cases}
\frac1{r^2},&\text{if the group \ $\Iso(X)$ is infinite;}\\
\ms_r(\mathcal D_{2n}\times \II,r),&\text{if \ $\Iso(X)$ is isomorphic to the dihedral group $\mathcal D_{2n}$ for some $n\ge 2$};\\
0,&\text{otherwise}.
\end{cases}
$$
\end{theorem}

Theorem~\ref{t5.3} suggests the following

\begin{problem}\label{7.4} Classify maximal symmetry numbers $\ms_r(X)$ of compact subsets $X$ of positive Lebesgue measure in $\IR^n$ for $n>2$. In particular, calculate
$\ms_r(X)$ for all regular convex polyhedra $X$ in $\IR^n$ for $n>2$. Calculate
$\ms_2(T)$ (or equivalently, $\ms_2(S_4\times \II)$) for the regular tetrahedron
$T$ in $\IR^3$.
\end{problem}

Observe that the last statement of Theorem~\ref{t5.2}, Theorem~\ref{t4.1}, and known
results on classification of regular polyhedra \cite[12.6]{Ber}
imply the following
asymptotic estimation first noticed by Ya.Vorobets.

\begin{proposition}\label{7.5} If $m\ge 2$ and $r=2^k$ for some $k\in\IN$, then
$$
\frac1{r^2}\ge \ms_r(X,r)\ge\frac1{r^2}\frac{2^m-r}{2^m-1}.
$$
for every regular convex polyhedron of dimension $\dim X\ge 2m-1$.
\end{proposition}

Maximal symmetry numbers $\ms_2(X,\mathcal S)$ of the finite sets $X$ of vertices of regular polyhedra in $\IR^3$ endowed with the counting measure $\mu:A\mapsto |A|$ and various families $\mathcal S$ of admissible symmetries were calculated by I.~Kluka in \cite{Kluka} and are presented in the following table, in which
\begin{itemize}
\item $\mathcal S^\pm$ is the family of  non-identity isometries of $X$;
\item $\mathcal S^+$ is the family of  non-identity orientation-preserving isometries of $X$;
\item $\mathcal S^-$ is the family of orientation-reversing isometries of $X$;
\item $\mathcal S^\pm_2$ is the family of non-identity involutive isometries of $X$;
\item $\mathcal S^+_2$ is the family of non-identity involutive orientation-preserving isometries of $X$;
\item $\mathcal S^-_2$ is the family of involutive orientation-reversing isometries of $X$.
\end{itemize}
\bigskip

\begin{center}
\begin{tabular}{| c| c| c |  c |  c | c |}
    \hline
$X$& tetrahedron&octahedron&cube&icosahedron&dodecahedron\\
\hline
$|X|$ &4&6&8&12&20\\
\hline
$\ms_2(X,\mathcal S^\pm)$ &$\mathsf 2$&$\mathsf 3$&$\mathsf 4$&$\mathsf 5$&$\mathsf 7$\\
$\ms_2(X,\mathcal S^\pm_2)$ &$\mathsf 2$&$\mathsf 3$&$\mathsf 4$&$\mathsf 5$&$\mathsf 7$\\
\hline
$\ms_2(X,\mathcal S^-)$ &$\mathsf 2$&$\mathsf 3$&$\mathsf 3$&$\mathsf 4$&$\mathsf 6$\\
$\ms_2(X,\mathcal S^-_2)$ &$\mathsf 2$&$\mathsf 3$&$\mathsf 3$&$\mathsf 5$&$\mathsf 6$\\
\hline
$\ms_2(X,\mathcal S^+)$ &$\mathsf 2$&$\mathsf 3$&$\mathsf 3$&$\mathsf 4$&$\mathsf 6$\\
$\ms_2(X,\mathcal S^+_2)$ &$\mathsf 2$&$\mathsf 2$&$\mathsf 2$&$\mathsf 4$&$\mathsf 6$\\
\hline
\end{tabular}
\end{center}
\bigskip

Observe that $\ms_2(X,\mathcal S)\ge\frac{|X|}4$ for all families of admissible symmetries $\mathcal S$ from this table. In particular, $\ms_2(X,\mathcal S^+_2)\ge\frac{|X|}4$ for every regular polyhedron $X$ in $\IR^3$. We do not know if the same inequality holds for the spheres $S^n$.

\begin{question} Is $\ms_r(S^n,\mathcal S^+_2)=\frac1{r^2}$ for
$n,r>1$? More precisely, does for every measurable $r$-coloring of the
$n$-dimensional sphere $S^n$ there exist a monochromatic subset of measure
$\ge\frac1{r^2}$, symmetric with respect to some line passing through the
center of the sphere?
\end{question}

It can be shown that $\ms_r(S^2,\mathcal S^+_2)\ge\frac1{2r^2}$ for every
$r\in\IN$ (use arguments analogous to those from \cite{Ba01}). Moreover,
 Theorems~\ref{t5.2}(2) implies that $$\ms_2(S^2,\mathcal S^+_2)\ge
\ms_2(S^2,\mathcal V)=\frac16\mbox{ \ and \ }\ms_3(S^2,\mathcal S_2^+)\ge \ms_3(S^2,\mathcal
V)\ge\frac1{18},$$ where $\mathcal V$ is the 3-element set consisting of
rotations on 180$^\circ$ around 3 pairwise orthogonal axes passing
through the center of the sphere $S^2$.

\section{The maximal symmetry numbers $\ms_r([0,1])$ of the interval $[0,1]$}

In this section we consider one of the most difficult and challenging open
problems related to our subject --- evaluating the maximal symmetry numbers
$\ms_r(\II)$ for the unit interval $\II=[0,1]$ in the real line $\IR$ endowed with the Lebesgue measure $\mu$ and the family $\mathcal S=\{s_c:c\in\IR\}$ of central symmetries $s_c:x\mapsto 2c-x$.

The number $\ms_r(\II)$ is tightly connected with the maximal symmetry numbers $\ms_r(n)$ of finite ordinals $n=\{0,1,\dots,n-1\}$ considered as subsets of integers $\IZ$ endowed with the counting measure $\mu:A\mapsto|A|$ and the family $\tilde S=\{\tilde s_c:c\in\IZ\}$ of admissible symmetries $\tilde s_c:x\mapsto c-x$. So, $\ms_r(n)$ is equal to the  maximal number $k\in\IN$ such that
for every $r$-coloring of the ordinal $n=\{0,\dots,n-1\}$ there
exists an $\tS$-symmetric monochromatic subset $A\subset n$ of cardinality $|A|\ge k$.

Note that the function $\ms_r(n)$ resembles the well known Van der
Waerden function $\mathsf{W}_r(n)$ assigning to each $n$ and $r$ the maximal $k$
such that for every $r$-coloring of $n=\{0,\dots,n-1\}$ there exists a monochromatic
arithmetic progression of length $k$. Since each arithmetic progression is
a symmetric set, we get $\mathsf{W}_r(n)\le \ms_r(n)$ for any $n,r\in\IN$. It is
know that for a given $r\ge 2$ the Van der Waerden function $\mathsf{W}_r(n)$ very
slowly approaches the infinity as $n\to \infty$. In the contract, the
function $\ms_r(n)$ grows linearly on $n$, that is, for
every $r\in\IN$ there exists a positive limit
$\lim\limits_{n\to\infty}\ms_r(n)/{n}$. It turns out that this limit is equal to the maximal symmetry number  $\ms_r(\II)$ of the unit interval $\II=[0,1]$. Let us recall that $\ms_r(\II)$ is equal to
the least upper bound of $\e>0$ such that for every measurable
$r$-coloring of $\II$ there exists a monochromatic symmetric subset
$A\subset\II$ of Lebesgue measure $\lambda(A)\ge\e$.

The exact values of the maximal symmetry numbers $\ms_r(\II)$ are not known. We merely have some lower and upper bounds for these numbers. Lower bounds for $\ms_r(\II)$ can be derived from known lower bounds for the function $\Delta:[0,1]\to[0,1]$ assigning to each $\e\in[0,1]$ the number
$$\Delta(\e)=\inf\{\ms(A):\mbox{$A\subset\II$ is a measurable subset of Lebesgue measure $\lambda(A)\ge\e$}\}.$$
The function $\Delta(\e)$ was introduced and studied by Martin and O'Bryant in \cite{MOB07} and \cite{MOB09}. In particular, they proved that $$\Delta(\e)\ge \mathsf S\cdot\e^2\mbox{ for every $\e\in[0,1]$}$$where $$\mathsf S=\inf_{f\in L_1\kern-1pt(\II)}\frac{\|f*f\|_\infty}{\|f\|^2_1}=\inf_{f\in L_1\kern-1pt(\II)}\frac{\sup_{x\in\IR}|\int_0^1f(t)f(x-t)dt|}{\big(\int_0^1|f(t)|dt\big)^2}$$and $L_1(\II)$ denote the space of measurable functions $f:\IR\to\IR$ with $\supp(f)=\{x\in\IR:f(x)\ne0\}\subset[0,1]$ and $\int_0^1|f(t)|dt<\infty$.

The exact value of the constant $\mathsf S$ is not known. Many authors proved lower and upper bounds for $\mathsf S$, see \cite{MOB07}, \cite{MOB09}, \cite{MV10}. At the moment the best bounds $$0.6374\le\mathsf S\le 0.7549$$ are due to Matolsci and Venuesa \cite{MV10}. The lower bound for the number $\mathsf S$ implies the following lower bound for the numbers $\ms_r(\II)$, which improves the earlier lower bound $\ms_r(\II)\ge1/(r^2+r\sqrt{r^2-r})$ proved in \cite{BVV00i} and \cite{BVV00e}.

\begin{theorem} $\ms_r(\II)\ge \Delta(\frac1r)\ge\mathsf S/r^2\ge 0.6374/r^2$ for all $r\in\IN$.
\end{theorem}

Upper bounds for the maximal symmetry numbers $\ms_r(\II)$ are obtained by the technique of blurred colorings, created in \cite{BVV00i} and \cite{BVV00e}.
To define such colorings we first look at ordinary colorings from a bit different
point of view.

Observe that every (measurable) $r$-coloring $\chi:X\to r$ of a space $X$ can be
identified with the collection $\{\chi_i\}_{i\in r}$ of (measurable) functions
$\chi_i:X\to\{0,1\}$ such that $\sum_{i\in r}\chi_i\equiv1$ and $\chi_i(x)=1$
iff $\chi(x)=i$ for $x\in X$. By a ({\it measurable}) {\it blurred
$r$-coloring} of $X$ we understand a collection $\{\chi_i\}_{i\in r}$ of
(measurable) functions $\chi_i:X\to[0,1]$ such that
$\sum_{i=1}^r\chi_i\equiv 1$.

Now assume that $\mathbf X=(X,\mu,\mathcal S)$ is a symmetry measure space whose admissible isometries $s\in\mathcal S$ are involutive. In this case the maximal symmetry number $\ms_r(\mathbf X)$ is equal to
$$
\inf_{\chi}\sup_{i\in r}\sup_{s\in\mathcal S}\int\limits_X\chi_i(x)\chi_i(s(x))d\mu, $$
where the infimum is taken over all measurable $r$-colorings
$\chi=\{\chi_i\}_{i\in r}$ of $X$. Taking this infimum over all measurable
blurred $r$-colorings of $X$, we obtain the definition of the blurred
maximal symmetry number $\bms_r(\mathbf X)$. Clearly, $\bms_r(\mathbf X)\le \ms_r(\mathbf X)$.

In particular, we can consider the blurred maximal symmetry numbers $\bms_r(n)$ and
$\bms_r([0,1])$ of discrete and continuous intervals (on the discrete
interval $n$ we consider the counting measure $\mu:A\mapsto|A|$).
The following theorem proved in \cite[\S6]{BVV00i} and \cite{BVV00e} describes
interplay between the blurred and ordinary maximal symmetry numbers.

\begin{theorem}\label{t6.1} Let $r\in\IN$. Then
\begin{enumerate}
\item $\ms_r(\II)=\lim\limits_{n\to\infty} \ms_r(n)/n=\inf\limits_{n\in\IN}\ms_r(n)/n$;
\item $\bms_r(\II)=\lim\limits_{n\to\infty}\bms_r(n)/n=\inf\limits_{n\in\IN}\bms_r(n)/n$;
\item $\bms_r(n)\le \ms_r(n)$ for every $n\in\IN$;
\item $\bms_r(\II)=\ms_r(\II)$ if $r>1$;
\item $\mathsf S/r^2\le \ms_r(\II)<1/r^2$ for $r>1$;
\item $\lim\limits_{r\to\infty}r^2\cdot\ms_r(\II)=\inf\limits_{r\in\IN}
r^2\cdot\ms_r(\II)=\mathsf S'$ for some constant $\mathsf S'$ such that $\mathsf S\le\mathsf S'<\frac56$;
\item $0.15935\le \ms_2(\II)<\frac5{24}=0.20833...$
\end{enumerate}
\end{theorem}

For small $n$ the values $\ms_2(n)$ can be easily found by routine computer calculations. On the other hand, the values $\bms_2(n)$ are not known even for small $n>3$.
For $n=4$ we know that $\bms_2(4)\le \frac5{24}$ and this explains why
the constant $\frac5{24}$ had appeared in Theorem~\ref{t6.1}. Answers (even
particular) to the following problems would be very interesting.

\begin{problem}\label{8.2} Calculate the values $\bms_r(n)$ (at least for
small $n$ and $r=2$).
\end{problem}

\begin{problem}\label{8.3} Calculate the values $\ms_r(\II)$ (at least for
$r=2$).
\end{problem}

\begin{problem}\label{8.4} Calculate or evaluate the constant $\mathsf S'$ from Theorem~\ref{t6.1}(6). Is $\mathsf S'>\mathsf S$?
\end{problem}

\begin{problem}\label{8.5} How quickly do the sequences
$\{\ms_r(n)/n\}_{n=1}^\infty$ and $\{\bms_r(n)/n\}_{n=1}^\infty$ converge to
their limit $\ms_r(\II)=\bms_r(\II)$?
\end{problem}

\section{The maximal symmetry numbers $\ms_r(G)$ of compact topological groups}

In this section we shall study the maximal symmetry numbers of symmetry measure spaces having Algebraic nature. Namely, we shall study the maximal symmetry numbers $\ms_r(G)$ of
compact topological groups $G$.

It is well-known that each compact topological group $G$ possesses a unique translation-invariant Borel probability measure $\mu:\dom(\mu)\to[0,1]$ called the {\em Haar measure}. This measure determines the $\sigma$-algebra $\dom(\mu)$ of $\mu$-measurable subsets, which includes all Borel subsets of $G$. Subsets that belong to $\dom(\mu)$ will be called {\em measurable}.

Each group $G$ will be endowed with the family of admissible symmetries $$\mathcal S=\{s_c:c\in G\}$$ containing involutive transformations $$s_c:x\mapsto cx^{-1}c,$$ called {\em central symmetries} of the group $G$.

For an Abelian group $G$ we shall also consider the wider family
$$\tS=\{\tilde s_c:c\in G\}$$of symmetries
$$\tilde s_c:x\mapsto c-x$$
called {\em quasi-central symmetries} of $G$.

If the Abelian group $G$ is 2-divisible (which means that the equation $x+x=a$ has a solution in $G$ for each $a\in G$), then $\tS=\mathcal S$.

For an Abelian compact topological group $G$ by $\ms_r(G)$ and $\tms_r(G)$ we shall denote the maximal symmetry numbers of the symmetry spaces $(G,\mu,\mathcal S)$ and $(G,\mu,\tS)$. These numbers take their values in the interval $[0,1]$.

Since each central symmetry $s_c$, $c\in G$, is involutive, for any measurable subset $A\subset G$ the set $A\cap s_c(A)$ is the largest measurable $\{s_c\}$-symmetric subset of $A$. The regularity of the Haar measure $\mu$ implies that the measure $\mu(A\cap s_c(A))$ continuously depends on $c$, so
$$\ms(A)=\max\{\mu(A\cap s_c(A)):c\in G\}$$is attainable by the compactness of $G$.

This implies that {\em for any measurable $r$-coloring of a compact topological group $G$ there is a monochromatic measurable $\mathcal S$-symmetric subset $A\subset G$ of Haar measure $\mu(A)\ge\ms_r(G)$ and $\ms_r(G)$ is the largest real number with that property}. A similar result holds for the maximal symmetry number $\tms_r(G)$.

The maximal symmetry numbers $\ms_r(G)$ and $\tms_r(G)$ of compact Abelian groups $G$ were evaluated in \cite{BVV00i}. Below by $G_2=\{x\in G:x+x=0\}$ we denote the closed subset of element of order 2 in $G$.

\begin{theorem}\label{t7.1} Let $r\ge 2$ be integer and $G$ be a compact Abelian group with Haar measure $\mu$.  \begin{enumerate}
\item $\frac1{r^2}\le \ms_r(G)\le\tms_r(G)$.
\item If $G$ is finite, then $$\frac1{r^2}<\tms_r(G)\le\frac1{r^2}+\Big(\frac1r-\frac1{r^2}\Big)\mu(G_2)+3\sqrt{2\ln(2r|G|)/|G|};$$
\item If $G$ is infinite, then $\tms_r(G)\le\frac1{r^2}+\big(\frac1r-\frac1{r^2}\big)\mu(G_2)$.
\end{enumerate}
\end{theorem}

The lower bounds in this theorem can be derived from the following proposition proved by the method of Harmonic Analysis in \cite[\S2]{BVV00i}.

\begin{proposition} Any measurable subset $A$ of a compact Abelian group $G$ contains a measurable $\mathcal S$-symmetric subset $B\subset A$ of Haar measure $\mu(B)\ge\mu(A)^2$. If $0<\mu(A)<1$, then $A$ contains a measurable $\tS$-symmetric subset $B\subset A$ of Haar measure $\mu(B)>\mu(A)^2$.
\end{proposition}

Theorem \ref{t7.1} can be completed by the following theorem of Yu.Zelenyuk \cite{Ze11s}.

\begin{theorem} Let $G$ be a finite Abelian group.
\begin{enumerate}
\item If $G$ contains a subgroup isomorphic to $\IZ_4^n$ for some $n\in\IN$, then $\ms_r(G)=\frac1{r^2}$ for $r=2^n$.
\item If $G$ contains no element of order 4, then $\ms_2(G)>\frac14$.
\end{enumerate}
\end{theorem}

This theorem suggests the following open:

\begin{problem} Given an integer number $r$, detect Abelian groups $G$ with $\ms_r(G)=\frac1{r^2}$.
\end{problem}

The evaluation of the maximal symmetry numbers $\ms_r(G)$ of non-commutative compact topological groups is  much more difficult and still remains an unfinished task.
We shall give some partial answer for finite and connected groups.

We start with the case of finite group $G$. To estimate the number $\ms_r(G)$, consider the following numbers:
$$
k(G)=\min_{g\in G}|\sqrt{g^2}|,\quad k_0(G)=|\sqrt 1|,\quad\text{and}\quad
k_1(G)=\min_{g\in 2G\setminus\{1\}}|\sqrt g|,
$$
where 1 is the neutral element of the group $G$, $2G=\{g^2 :g\in
G\}\subset G$ and $\sqrt g=\{x\in G:x^2=g\}$ for $g\in G$. Let also $m(G)$
be the maximal cardinality of a subgroup of $G$, contained in the set $2G$.

\begin{theorem}[\cite{Ba99v}]\label{t7.5} For every finite group $G$ and every
$r\in\IN$ we have the estimates:
\begin{enumerate}
\item
$\ms_r(G)\le\frac1{r^2}+(\frac1r-\frac1{r^2})\frac{k_0(G)}{|G|}+3\sqrt{2\ln(2r|G|)/|G|}$;
\item $\ms_r(G)\ge \frac{k(G)m(G)}{|G|}\cdot\frac1{r^2}$;
\item $\ms_r(G)\ge
\frac{k_1(G)m(G)}{|G|}\frac1{r^2}+\frac{k_0(G)-k_1(G)}{|G|}\frac1{r}$;
\item $\ms_r(G)\ge\frac1{r^2}$ of the group $G$ has odd cardinality $|G|$.
\end{enumerate}
\end{theorem}

Similar lower and upper bound were also obtained by M.Korostenski and Yu.Zelenyuk in \cite{KZ10}.

For some groups (for example, dihedral groups) the third item of Theorem~\ref{t7.5}
gives a better lower bound than the second one.
For any finite Abelian group $G$ we have $m(G)=|2G|$ and
$k(G)=|B(G)|=\frac{|G|}{m(G)}$ and thus
$\ms_r(G)\ge\frac1{r^2}$ (which agrees with Theorem~\ref{t7.1}). The last statement of Theorem~\ref{t7.5} follows from the third one and the observation (due to I.V.Protasov) that each finite group $G$ of odd order $|G|$ has $2G=G$, $m(G)=|G|$ and $k(G)=k_0(G)=k_1(G)=1$. Yet, for
Abelian groups of odd order we have the strict lower bound that follows from Theorem~\ref{t7.1}(2):

\begin{corollary} Each finite Abelian group $G$ of odd cardinality has the maximal symmetry numbers $\ms_r(G)=\tms_r(G)>\frac1{r^2}$ for all finite $r\ge 2$.
\end{corollary}

The proof of this fact exploits the Fourrier
transformation, which does not work properly in the
non-Abelian case. So, the following problem appears naturally:

\begin{question} Is $\ms_r(G)>\frac1{r^2}$ for every group $G$ of
odd order and every $r>1$?
\end{question}

Surprisingly, but for non-commutative groups of even order the inequality
$\ms_r(G)\ge\frac1{r^2}$ does not hold anymore: it was noticed in \cite{Ba99v}
that $\ms_r(Q_8)=\frac18<\frac1{2^2}$, where $Q_8=\{\pm1,\pm i,\pm j,\pm
k\}$ is the group of quaternions. Nonetheless the following question
(related to Theorem~\ref{t8.3}) remains open.

\begin{question}[\cite{Ba99v}] Does for every $\e>0$ and $r>1$ there
exist a finite group $G$ with $\ms_r(G)<\e$?
\end{question}

Next, we  evaluate the maximal symmetry numbers $\ms_r(G)$ of
non-commutative connected compact topological groups. We begin with connected Lie
groups.

\begin{theorem}[\cite{Ba01}]\label{t7.9} If $r\in\IN$ and $G$ is a compact
connected Lie group, then
$$
\frac1{r^2}\cdot \frac{k(G)}{2^{\dim(G)}}\le \ms_r(G)\le\frac1{r^2},
\quad\text{where $k(G)=\min_{g\in G}|\sqrt{g^2}|$}.
$$
\end{theorem}

It is interesting to note that $k(\IT^n)=2^n$ for the $n$-th power of the circle $\IT=\{z\in\IC:|z|=1\}$ and thus Theorem 5.2 gives
the exact value for $\ms_r(\IT^n)$. For arbitrary
connected compact group we have a weaker lower bound.

\begin{theorem}[\cite{Ba01}] If $r\in\IN$ and $G$ is a
non-degenerate connected compact group, then
$$
\frac1{r^2}\cdot\frac1{2^{\dim G}}\le \ms_r(G)\le\frac1{r^2}.
$$
\end{theorem}

\begin{question}[\cite{Ba01}] Is there a connected compact group
$G$ with $\ms_r(G)<\frac1{r^2}$ for some $r\in\IN$? In particular, what
is the value of $\ms_2(SO(3))$? What is the value of $\ms_2(S^3)$ where $S^3$ is the multiplicative group of quaternions with unit norm?
\end{question}

Note that Theorems~\ref{t7.9} yields the estimate $\frac1{16}\le
\ms_2(S^3)\le\frac14$.

\section{The cardinal symmetry numbers $\MS_r^+(G)$ of groups}

In this section we consider symmetry measure spaces $\mathbf X=(X,\card,\mathcal S)$ endowed with the cardinal-valued measure $\card:\mathcal P(X)\to[0,|X|]$ defined on the $\sigma$-algebra $\mathcal P(X)$ of all subsets of $X$ and assigning to each subset $A\subset X$ its cardinality $\card(A)=|A|$. A bit more informative is the measure $\card^+:\mathcal P(X)\to[1,|X|^+]$ assigning to $A\subset X$ the successor $|A|^+$ of the cardinality $|A|$ of $A$.

We shall be interested in the maximal symmetry numbers $\ms_r(G,\card,\mathcal S)$ and $\ms_r(G,\card^+,\mathcal S)$ of groups $G$ endowed with the family $\mathcal S=\{s_c:c\in G\}$ of central symmetries $s_c:x\mapsto cx^{-1}c$. These symmetry Ramsey numbers will be denoted by $\MS_r(G)$ and $\MS^+_r(G)$, respectively, and will be referred to as {\em cardinal symmetry numbers} of $G$. It is easy to see that
$$\MS_r(G)\le\MS^+_r(G)\le(\MS_r(G))^+.$$

It follows from the definition that $\MS_r^+(G)$ is equal to the maximal cardinal $\kappa$ such that for each cardinal $\lambda<\kappa$ any each $r$-coloring of $G$ there is a monochromatic $\mathcal S$-symmetric subset $A\subset G$ of cardinality $|A|>\lambda$.

For each finite group $G$ we get
$$\MS_r^+(G)=1+\MS_r(G)=1+|G|\cdot\ms_r(G).$$
So, now we shall concentrate on infinite groups $G$. We start with the following result of Protasov \cite{Pr96}.

\begin{theorem}\label{t8.1} For any infinite Abelian group $G$ and a finite cardinal $r$ the cardinal symmetry numbers $\MS_r(G)$ and $\MS_r^+(G)$ are infinite.
\end{theorem}

This theorem follows from its more precise version taken from \cite{BP01}:

\begin{theorem}\label{t8.2} If $r,n$ are natural numbers, then for every
$r$-coloring of an Abelian group $G$ of cardinality $|G|\ge r^2n$ there exists an $\mathcal S$-symmetric monochromatic subset $A\subset G$ of cardinality $|A|\ge n$.
\end{theorem}

A similar result holds also for finite non-commutative groups.
Below $\pi(n)$ is the function assigning to each
$n$ the quantity of prime numbers not exceeding $n$.

\begin{theorem}\label{t8.3} If $r,n$ are natural numbers, then for every
$r$-coloring of a finite group $G$ with $|G|>r^2(n!)^{\pi(n)}$
there exists a symmetric monochromatic subset $A\subset G$ of cardinality
$|A|>n$.
\end{theorem}

This theorem follows from Theorem~\ref{t8.2} and the known fact stating that each
finite group $G$ with $|G|>(n!)^{\pi(n)}$ contains an Abelian
subgroup $H$ of cardinality $|H|>n$. In its turn, this fact can be easily
derived from the Silov Theorem 11.1 and Theorem 16.2.6 of \cite{KM}. Note
that in contrast to Theorem~\ref{t8.2}, Theorem~\ref{t8.3} says nothing about infinite
groups.

\begin{problem}[Protasov]\label{pr8.4} Are the cardinals $\MS_r(G)$ and $\MS^+_r(G)$ infinite for every infinite group $G$ and every finite cardinal $r$? Equivalently, is it true that for every finite coloring of an infinite group $G$ there exist a monochromatic symmetric subset $A\subset G$ of arbitrary large finite cardinality?
\end{problem}

Note that this problem does not reduce to the Abelian case since there
exist infinite periodic groups containing no big Abelian subgroups (for
example, for every odd $d\ge 665$ and every $m\ge 2$ the free Bernside
group $B(m,d)$ is infinite and contains no Abelian subgroup of cardinality $>d$,
see \cite{Ad}). In the meantime, by the theorem of Kargapolov-Hall-Kulatilani
(see \cite{NRR}), each locally finite infinite group does contain an
infinite Abelian subgroup and thus for locally finite groups Problem~\ref{pr8.4}
has affirmative answer.

Next, given an Abelian group $G$, we detect numbers $r$ for which the cardinal symmetry number $\MS_r^+(G)$ is uncountable. Observe that the cardinal $\MS_r^+(G)$ is uncountable if and only if for each $r$-coloring of $G$ there is an infinite monochromatic symmetric subset of $G$.

To describe numbers $r$ with uncountable $\MS_r(G)$ we need to recall the definition of $p$-ranks $r_p(G)$ of an Abelian group. For a non-negative integer $p\ge 0$ by the {\em $p$-rank} of an abelian group $G$ we understand the largest cardinal number $\kappa$ such that the direct sum $\oplus_{\alpha\in\kappa}\IZ_p$ of $\kappa$ many copies of the cyclic group $\IZ_p$ is isomorphic to a subgroup of the group $G$. For $p=0$ we assume that the group $\IZ_0=\IZ$ is infinite cyclic.

\begin{theorem}[\cite{BP01}]\label{t8.5} For an Abelian group $G$ and a cardinal number $r$ the cardinal $\MS_r^+(G)$ is uncountable if and only if one of the following conditions is satisfied:
\begin{enumerate}
\item $r\le r_0(G)$ and the group $G$ is finitely-generated;
\item $r\le r_0(G)+1$ and the group $G$ is countable, infinite-generated, and has finite 2-rank $r_2(G)$;
\item $r\le\max\{\log|G|,r_2(G)\}$ and the group $G$ is uncountable or has infinite 2-rank $r_2(G)$.
\end{enumerate}
\end{theorem}

In this theorem for a cardinal number $\kappa$ we put $\log(\kappa)=\min\{\lambda:\kappa\le 2^\lambda\}$.

Since $\MS_r^+(G)\le|G|^+$, Theorems~\ref{t8.1} and \ref{t8.5} yield the complete description of the cardinal function $\MS_r^+(G)$ for countable Abelian groups $G$. For uncountable Abelian groups $G$ there is a non-trivial upper bound for $\MS_r^+(G)$, proved by I.Protasov \cite{Pr99}:

\begin{theorem}\label{t8.6} For any uncountable Abelian group $G$ with $r_2(G)<|G|$ and every cardinal $r\ge 2$ we have $\MS_r^+(G)\le\MS_2^+(G)\le|G|$, which follows from the fact that the group $G$ admits a 2-coloring of $G$ without monochromatic $\mathcal S$-symmetric subsets of cardinality $|G|$.
\end{theorem}

This theorem is not true for non-commutative groups because of the following example due to Yu.Gryshko and A.Khelif \cite{GK05}.

\begin{example} For any infinite field $K$ of characteristic $\ne 2$, the special linear group $G=SL(K)$ has cardinal symmetry number $\MS_2^+(G)=|G|^+$, which means that for any 2-coloring of the group $G=SL(K)$ there is a monochromatic symmetric subset of cardinality $|G|$.
\end{example}

\begin{problem} Let $G=SL(K)$ for an infinite field $K$. Is $\MS_r^+(G)=|G|^+$ for all finite cardinals $r$?
\end{problem}

Uncountable groups $G$ with $\MS_2^+(G)\le |G|$ were characterized by Yu.Gryshko and A.Khelif in \cite{GK05}.

\begin{theorem} An uncountable group $G$ has $\MS_2^+(G)\le|G|$ if and only if the subgroup $L(G)=\{(x,y)\in G\times G:xy\in[G,G]\}$ can be written as the union $L(G)=\cup\mathcal H$ of a linearly ordered chain $\mathcal H$ of subgroups such that for any subgroup $H\in\mathcal H$ we get $|H|<|G|$ and  $|\{x\in G:\exists (u,v)\in H\;x^2=uv\}|<|G|$. If the cardinal $|G|$ is regular, then the latter condition is equivalent to $|\{x\in G:x^2=a\}|<|G|$ for all $a\in[G,G]$.
\end{theorem}

Lower estimates for the cardinal symmetry numbers $\MS_r(G)$ of uncountable groups $G$ can be proved using Erd\"os-Rado partition relations for ordinals.

Let us recall that for ordinals $\kappa$, $\lambda$, a cardinal $r$, and a natural number $n$, the symbol
 $\kappa\to(\lambda)^n_r$ means that for any $r$-coloring of the set $[k]^n$ of $n$-element subsets of $\kappa$ there is a subset $\Lambda\subset \kappa$ of order type $\lambda$ such that the subset $[\Lambda]^n$ is monochrome, see \cite{EHMR}.

The following theorem was announced by T.Banakh and I.Protasov in \cite{BP01} and proved by Yu.Gryshko and A.Khelif in \cite{GK05}.

\begin{theorem}\label{t8.10} If $\kappa,\lambda,r$ are cardinals such that $\kappa\to(\lambda+1)^2_r$ then any uncountable Abelian group $G$ of cardinality $|G|\ge\kappa$ has cardinal symmetry number $\MS^+_r(G)>\lambda$, which means that for any $r$-coloring of $G$ there is a monochrome symmetric subset of cardinality $\ge\lambda$.
\end{theorem}

It follows from the famous Erd\"os-Rado stepping up lemma \cite{EHMR} that $(2^{<\kappa})^+\to (\kappa+1)^2_r$ for every cardinals $r<\kappa\ge\aleph_0$, where $2^{<\kappa}=\sup\{2^\tau:2^\tau:\tau<\kappa\}$. This observation and Theorem~\ref{t8.10} imply

\begin{corollary}\label{c8.11} For every infinite cardinal $\kappa$, and a cardinal $r<\kappa$, any Abelian group $G$ of cardinality $|G|>2^{<\kappa}$ has cardinal symmetry number $\MS^+_r(G)>\kappa$.
\end{corollary}

Observe that GCH, the Generalized Continuum Hypothesis, is equivalent to the
assumption that the equality $\kappa=2^{<\kappa}$ holds for all infinite cardinals $\kappa$.
Then Theorem~\ref{t8.6} and Corollary~\ref{c8.11} imply the following GCH-theorem:

\begin{theorem} Under GCH, for any uncountable Abelian group $G$ we get
$$\MS_r^+(G)=\begin{cases}
|G|^+&\mbox{if $r<|G|=r_2(G)$,}\\
|G|&\mbox{if $r<|G|$ and $r_2(G)<|G|$,}\\
1&\mbox{if $r\ge|G|$.}
\end{cases}
$$
\end{theorem}

We do not know if this result is true in ZFC. The following theorem of I.Protasov suggests a positive
answer to this question.

\begin{proposition} For any uncountable Abelian group $G$ with $r_2(G)<|G|$ the equality $\MS^+_r(G)=|G|$ holds for every $r\in\{2,3\}$.
\end{proposition}

This proposition  can be easily derived from the next modification of
Lemma 1 of \cite{Pr99}.

\begin{lemma} Let $G$ be an infinite Abelian group, $k<|G|$, and
$H\subset G$ be a subgroup with $|H|\le k$. For every $k$-coloring
$\chi:G\to k$ of the group $G$ one of the following statements holds:
\begin{enumerate}
\item $G$ contains a monochromatic symmetric subset of cardinality $\ge k$;
\item there exists $g\in G$ such that
$\chi(2H+g)\cap\chi(2H-g)=\emptyset$.
\end{enumerate}
\end{lemma}

In contrast to the Abelian case, the problem of calculation of the cardinal symmetry numbers $\MS^+_r(G)$ of  non-commutative groups is much less studied. Definitive answers (due to  Yu.Gryshko and A.Khelif) are known only for $r=2$.

\begin{theorem}[\cite{GK05}]\label{t8.15} If an infinite group $G$ has cardinal symmetry number $\MS_2^+(G)\le\aleph_0$, then $G$ is countable and is either locally finite or virtually cyclic.
\end{theorem}

Let us recall that a group $G$ is called
\begin{itemize}
\item {\em locally finite} if each finitely-generated subgroup of $G$ is finite;
\item {\em virtually cyclic} if $G$ contains a cyclic subgroup of finite index.
\end{itemize}

The paper \cite{GK05} contains also a characterization of locally finite and virtually cyclic groups $G$ with $\MS_2^+(G)\le\aleph_0$.

\begin{theorem}\label{t8.16} A countable locally finite group $G$ has cardinal symmetry number $\MS_2^+(G)\le\aleph_0$ if and only if for any element $a\in[G,G]$ in the commutator subgroup $[G,G]$ of $G$ the set $\sqrt{a}=\{x\in G:x^2=a\}$ is finite.
\end{theorem}

\begin{theorem}\label{t8.17} A virtually cyclic group $G$ with normal cyclic subgroup $H$ of finite index has cardinal symmetry number $\MS_2^+(G)\le\aleph_0$ if and only if either $G$ coincides with the centralizer $C_H=\{g\in G:\forall h\in H\;\;gh=hg\}$ of $H$ or else there is an element $g\in G$ whose square $g^2$ does not belong to the subgroup generated by the set $[G,G]\cup\{c^2:c\in C_H\}$.
\end{theorem}

Theorem~\ref{t8.15} implies that the free group $F_2$ with two generators has uncountable cardinal symmetry number $\MS_2^+(F_2)$, which means that for each 2-coloring of $F_2$ there is an infinite  monochromatic symmetric subset of $F_2$.

\begin{problem} Is $\MS_r^+(F_2)$ uncountable for any finite cardinal $r$?
\end{problem}

\section{The dimension symmetry numbers $\mds_r(X)$ of some metric spaces}

In this section we consider the maximal symmetry numbers of metric spaces $X$ endowed with the (non-additive) measure $\asdim:\mathcal P(X)\to \IZ\cup\{+\infty\}$ assigning to each subset $A\subset X$ its asymptotic dimension $\asdim(A)$.

We say that for an integer number $n\ge-1$, a subset $A$ of a metric space $(X,d)$ has asymptotic dimension $\asdim(A)\le n$ if for each finite real number $\e$ there is a cover $\U$ of $A$ with finite $\mesh(\U)=\sup\limits_{U\in\U}\diam(U)$ such that for some bounded subset $B\subset X$ and each point $x\in X\setminus B$ the $\e$-ball $B(x,\e)=\{y\in X:d(x,y)<\e\}$ centered at $x$ meets at most $n+1$ elements of the cover $\U$.

The {\em asymptotic dimension} $\asdim(A)$ of $A\subset X$ is equal to the smallest integer number $n\ge-1$ such that $\asdim(A)\le n$. If such number $n$ does not exist, then we put $\asdim(A)=\infty$. It follows that a subset $A\subset X$ is bounded if and only if $\asdim(A)=-1$.

It is clear that $\asdim(A)\le\asdim(B)$ for any subsets $A\subset B$ of $X$. Thus the asymptotic dimension can be considered as a (non-additive) measure $\asdim:\mathcal P(X)\to [-1,\infty]$.
This measure is idempotent in the sense that $\asdim(A\cup B)=\max\{\asdim(A),\asdim(B)\}$ for all subsets $A,B\subset X$, see \cite{Roe}.

Assume that a metric space $X$ is given with a family $\mathcal S$ of admissible symmetries that  preserve the asymptotic dimension of subsets of $X$. This happens if each admissible symmetry $s\in\mathcal S$ is a {\em quasi-isometry} in the sense that there is a finite positive constant $C$ such that
$$\frac1C d_X(x,x')-C\le d_Y(f(x),f(x'))\le C\cdot d_X(x,x')+C$$for all $x,x'\in X$.

In this case the triple $\mathbf X=(X,\asdim,\mathcal S)$ is a symmetry measure space whose maximal symmetry number $\ms_r(X,\asdim,\mathcal S)$ will be called the {\em dimension symmetry number} of $(X,\mathcal S)$ and will denoted by $\mds_r(X,\mathcal S)$ or $\mds_r(X)$ if the family $\mathcal S$ of admissible symmetries is clear from the context. In some cases (for example for the hyperbolic plane) it is necessary to restrict the domain of the definition of the measure $\asdim$ to the $\sigma$-algebra $\Borel$ of Borel subsets of $X$. In this case we obtain a symmetry measure space $(X,\asdim|\Borel,\mathcal S)$ whose maximal symmetry number $\ms_r(X,\asdim|\Borel,\mathcal S)$ will be denoted by $\mds_r^B(X,\mathcal S)$ or just $\mds_r^B(X)$ if the family $\mathcal S$ of admissible symmetries is clear from the context.

It is clear that $\mds_r(X)\le \mds_r^B(X)$ for any $r\in\IN$. If the metric space $X$ is discrete, then each subset of $X$ is Borel and hence $\mds_r(X)=\mds_r^B(X)$.

It follows from the definition of $\mds^B_r(X)$ that for each number $n\le\mds_r^B(X)$ and each Borel $r$-coloring of $X$ there is a monochromatic $\mathcal S$-symmetric subset $S\subset X$ of asymptotic dimension $\asdim(S)\ge n$, and $\mds_r^B(X)$ is the smallest number with that property.

Since a subset $A$ of a metric space $X$ is unbounded if and only if $\asdim(A)\ge 0$, the space $X$ has dimension symmetry number $\mds_r^B(X)\ge 0$ if and only if for any Borel $r$-coloring, the space $X$ contains an unbounded monochromatic $\mathcal S$-symmetric subset. If all bounded subsets in $X$ are finite, then $\mds_r^B(X)=\mds_r(X)\ge 0$ if and only if $\MS_r(X)$ is infinite.

First we discuss the dimension symmetry numbers $\mds_r(\IR^n)$ and $\mds^B_r(\IR^n)$ of the Euclidean spaces $\IR^n$ endowed with the family $\mathcal S=\{s_c:c\in\IR^n\}$ of central symmetries $s_c:x\mapsto 2c-x$.

The Euclidean space $\IR^n$ contains the discrete subspace $\IZ^n$, endowed with the induced Euclidean metric and the family $\mathcal S=\{s_c:c\in\IZ^n\}$ of central symmetries.

\begin{theorem}[\cite{BC?}]\label{t9.1} For any natural numbers $n,r$ we get
$$\mds_r(\IR^n)=\mds_r^B(\IR^n)=\mds_r(\IZ^n)=n-r+1.$$
This follows from the following facts:
\begin{enumerate}
\item for any $r$-coloring of $\IZ^n$ there is a monochromatic $\mathcal S$-symmetric subset of asymptotic dimension $\ge n-r+1$;
\item there is a Borel $r$-coloring of $\IR^n$ without monochromatic $\mathcal S$-symmetric subsets of asymptotic dimension $>n-r+1$.
\end{enumerate}
\end{theorem}

Next, we consider the dimension symmetry numbers $\mds_r(\IH^2)$ and $\mds^B_r(\IH^2)$ of the hyperbolic plane $\IH^2$ endowed with the family $\mathcal S=\{s_c:c\in\IH^2\}$ of central symmetries. For each point $c\in\IH_2$ the central symmetry $s_c:\IH^2\to\IH^2$ assigns to each point $x\in\IH^2$ the unique point $s_c(x)\in\IH^2$ such that $c$ is the midpoint of the segment $[x,s_c(x)]$.

The following unexpected result proved in \cite{BDR10} shows that by its Ramsey-theoretic properties, the hyperbolic plane substantially differs from its Euclidean counterpart.

\begin{theorem} For any $r\in\IN$ we get $\mds_r^B(\IH^2)\ge 0$, which means that for any Borel $r$-coloring of the hyperbolic plane $\IH^2$ there is an unbounded monochromatic symmetric subset $A\subset\IH^2$.
\end{theorem}

We do not know if this theorem remains true for arbitrary (not necessary Borel) colorings of $\IH^2$.

\begin{problem} Is $\mds_r(\IH^2)\ge 0$ for any $r\ge 2$?
\end{problem}

The answer to this problem is affirmative for $r=2$. This follows from the subsequent theorem that can be proved by analogy with Theorem 3.1 of \cite{BDR10}:

\begin{theorem}\label{t9.4} The family $\mathcal S=\{s_c:c\in\IH^2\}$ of central symmetries of the hyperbolic plane $\IH^2$ contains a 3-element subfamily $\mathcal S'$ such that
for the any 2-coloring of $\IH^2$ there is a monochromatic $\mathcal S'$-subset $A\subset \IH^2$ of asymptotic dimension $\asdim(A)=2$. Consequently, $\mds_2(\IH^2)=2$.
\end{theorem}

A bit weaker result holds also for the Euclidean plane $\IR^2$:

\begin{theorem}  The family $\mathcal S=\{s_c:c\in\IR^2\}$ of central symmetries of the Euclidean plane $\IR^2$ contains a 3-element subfamily $\mathcal S'$ such that for each 2-coloring of $\IR^2$ there is a monochromatic $\mathcal S'$-symmetric subset $A\subset\IR^2$ of asymptotic dimension $\asdim(A)\ge 1$.
\end{theorem}

This theorem is a particular case of the following result proved in \cite{BC?}:

\begin{theorem}\label{t9.6} For every finite numbers $r\le n$,  the family $\mathcal S=\{s_c:c\in\IZ^n\}$ of central symmetries of $\IZ^n$ contains a finite subfamily $\mathcal S'\subset\mathcal S$ of cardinality $|\mathcal S'|<2^n$ such that for any $r$-coloring of $\IZ^n$ there is a monochromatic $\mathcal S'$-symmetric subset $A\subset\IZ^n$ of asymptotic dimension $\asdim(A)\ge n-r+1$.
\end{theorem}

\section{Centerpole numbers $\cp_r^\lambda(X)$ of symmetry measure spaces}

Theorems~\ref{t9.4}--\ref{t9.6} motivate the following definition. Let $\mathbf X=(X,\mu,\mathcal S)$ be a symmetry measure space, $r$ be a cardinal, and $\lambda$ be an element of the linearly ordered set $[\mu(\emptyset),\mu(X)]$. A subfamily $\mathcal S'\subset\mathcal S$ will be called {\em $(r,\lambda)$-centerpole} if for any measurable $r$-coloring of $X$ there is a monochromatic $\mathcal S'$-symmetric subset $A\in\dom(\mu)$ of measure $\mu(A)\ge\lambda$.
The minimal cardinality $|\mathcal S'|$ of a $(r,\lambda)$-centerpole subset $\mathcal S'\subset S$ will be denoted by  $\cp_r^\lambda(\mathbf X)$ and called the {\em $(r,\lambda)$-centerpole number} of $\mathbf X$. If no $(r,\lambda)$-centerpole set exists (which happens if $\lambda>\ms_r(\mathbf X)$), then we put $\cp_k^\lambda(\mathbf X)=\infty$ and assume that $\infty$ is greater than any cardinal.

Observe that the symmetry Ramsey number $\ms_r(\mathbf X)$ can be expressed via the centerpole numbers by the formula:
$$\ms_r(\mathbf X)=\sup\{\kappa:\forall \lambda<\kappa\;\;\;\cp_r^{\lambda^+}(\mathbf X)<\infty\}.$$

Now we consider the following general problem:

\begin{problem} Given a symmetry measure space $\mathbf X=(X,\mu,\mathcal S)$ calculate its $(r,\lambda)$-centerpole numbers $\cp_r^\lambda(\mathbf X)$.
\end{problem}

We shall partly answer this question for Euclidean space $\IR^n$ and their discrete subspaces $\IZ^n$, viewed as symmetry measure spaces endowed with the measure $\asdim$ and the family $\mathcal S=\{s_c:c\in X\}$ of central symmetries. In fact, some information on the centerpole numbers $\cp_r^\lambda(X)$ is available only for $\lambda=0$.

Observe that for a symmetry measure space $\mathbf X=(X,\mu,\mathcal S)$ the centerpole number $\cp^0_r(X)$ is equal to the smallest cardinality $|\mathcal S'|$ of a subfamily $\mathcal S'\subset \mathcal S$ such that for each $r$-coloring of $X$ there is an unbounded monochromatic $\mathcal S'$-symmetric subset $A\subset X$. Here we remark that a subset $A\subset X$ has asymptotic dimension $\asdim(A)\ge0$ if and only if it is unbounded.

In the following theorem we present all available information on exact values of the centerpole numbers $\cp_r^0(\IR^n)$ and $\cp_r^0(\IZ^n)$ of the metric spaces $\IR^n$ and $\IZ^n$.

\begin{theorem}[\cite{BC11}]\label{t10.2}
\begin{enumerate}
\item $\cp_1^0(\IR^n)=\cp_1^0(\IZ^n)=1$ for any $n\ge 1$;
\item $\cp_2^0(\IR^n)=\cp_2^0(\IZ^n)=3$ for any $n\ge 2$;
\item $\cp_3^0(\IR^n)=\cp_3^0(\IZ^n)=6$ for any $n\ge 3$;
\item $8\le \cp_4^0(\IR^n)\le \cp_4^0(\IZ^n)\le 12$ for any $n\ge 4$;
\item $\cp_4^0(\IR^4)=\cp_4^0(\IZ^4)=12$;
\end{enumerate}
\end{theorem}

For larger numbers $r\le n$ we merely have some lower and upper bounds proved in \cite{BC11}:

\begin{theorem}[\cite{BC11}]\label{t10.3} For any $r\le n$ we have the following bounds:
\begin{enumerate}
\item $\cp_r^0(\IR^n)\le \cp^0_r(\IZ^n)\le \cp^0_r(\IZ^r)\le 2^r-1-\max\limits_{s\le r-2}\binom{r-1}{s-1}$;
\item $\cp^0_r(\IR^r)\ge \frac12(r^2+3r-4)$;
\item $\cp^0_r(\IR^n)\ge r+4$;
\item $\cp^0_r(\IR^n)<\cp^0_{r+1}(\IR^{n+1})$ and $\cp^0_r(\IZ^n)<\cp_r^0(\IZ^{n+1})$.
\end{enumerate}
\end{theorem}

The binomial coefficient $\binom{r}{s}$ in statement (1) is equal to $\frac{k!}{s!(n-s)!}$ if $s\in\{0,\dots,r\}$ and zero otherwise. By Theorem~\ref{t9.1}, $\mds_r(\IR^n)=\mds_r(\IZ^n)$. We do not know if the same is true for the centerpole numbers.

\begin{problem} Is $\cp_r^0(\IR^n)=\cp_r^0(\IZ^n)$ for all $r\le n$?
\end{problem}

Since each central symmetry $s_c:x\mapsto 2c-x$ of $\IR^n$ can be identified with its fixed point $c$ we can identify $k$-centerpole subsets with subsets of $\IR^n$ and study their geometric structure. We shall say that a subset $C\subset\IR^n$ is {\em $r$-centerpole} if the family of central symmetries $\mathcal S_C=\{s_c:c\in C\}$ is $(r,0)$-centerpole, which means that $\ms_r(\IR^n,\asdim,\mathcal S_C)\ge0$. The following theorem proved in \cite{Ba98} and \cite{BC11} describes the geometry of $r$-centerpole sets for small $r$.

\begin{theorem}\label{t10.5} \begin{enumerate}
\item  A subset $C\subset \IR^n$, $n\ge 1$, is 1-centerpole if and only if $C$ is not empty.
\item  A finite subset $C\subset\IR^n$, $n\ge 2$, is 2-centerpole if and only if $C$ contains a triangle, i.e., a 3-element subset $\{a,b,c\}$ that does not lie on a line.
\item A 6-element subset $C\subset \IR^3$ is 3-centerpole if and only if $C$ is an affine octahedron, i.e., $C=\{c+e_1,c-e_1,c+e_2,c-e_2,c+e-3,c-e_3\}$ for some point $c\in\IR^n$ and some linearly independent vectors $e_1,e_2,e_3$ in $\IR^n$.
\item There is a finite 3-centerpole subset of $\IR^3$ which does not contain any affine octahedron.
\item A subset $C\subset\IR^n$ is $(r+1)$-centerpole if it contains an affine copy of the $\binom{r}{s}$-sandwich $\Xi^{r}_{s}$ for some $s\le r-2$.
\end{enumerate}
\end{theorem}

The $\binom{r}{s}$-sandwich $\Xi^r_s$ mentioned in the preceding theorem is a special subset of $\IZ^{r+1}$, defined as follows.

Let $\bold 2=\{0,1\}$ and $\bold 2^r$ be the discrete $r$-cube in $\IZ^r\subset\IR^r$. For any real number $s$ consider the subsets
$$\bold 2^r_{<s}=\big\{(x_i)\in\bold 2^r:\sum_{i=1}^rx_i<s\big\}\mbox{ \ and \ }\bold 2^r_{>s}=\big\{(x_i)\in\bold 2^r:\sum_{i=1}^rx_i>s\big\}$$called the $s$-slices of the $r$-cube $\bold 2^r$. The subset
$$\Xi^r_s=\big(\{-1\}\times \bold 2^r_{<s}\big)\cup\big(\{0\}\times \bold 2^r_{<r}\big)\cup \big(\{1\}\times\bold 2^r_{>s}\big)$$of $\IZ\times\IZ^r$ is called the {\em $\binom{r}{s}$-sandwich}. In \cite{BC11} it is proved that the $\binom{r}{s}$-sandwich is a $(r+1)$-centerpole set in $\IZ^{r+1}$ for every $s\le r-2$. The minimal cardinality $$\min_{s\le r-2}|\Xi^r_s|=2^{r+1}-1-\max_{s\le r-2}{\textstyle\binom{r}{s}}$$ of such sandwiches yields the upper bound from Theorem~\ref{t10.3}(1).

Observe that
\begin{itemize}
\item $\Xi^0_{-2}=\{(1,0)\}$ is a singleton in $\IZ\times\IZ^0=\IZ\times\{0\}$;
\item $\Xi^1_{-1}=\{(0,1),(1,0),(1,1)\}$ is the unit square without a vertex in $\IZ^2$;
\item $\Xi^2_{0}=\{(0,0,0),(0,0,1),(0,1,0),(1,0,1),(1,1,0),(1,1,1)\}$ is the unit cube without two opposite vertices in $\IZ^3$;
\item $\Xi^3_0$ is the unit cube without two opposite vertices in $\IZ^4$;
\item $\Xi^3_1$ is a 12-element subset of $\IZ^4$ whose slices $\{-1\}\times \mathbf 2^3_{-1}$, $\{0\}\times \mathbf 2^3_{<3}$ and $\{1\}\times\mathbf 2^3_{>1}$ have 1, 7 and 4 points, respectively.
\end{itemize}

It is interesting to observe that the sandwich
\begin{itemize}
\item $\Xi^0_{-2}$ is affinely equivalent to a singleton in $\IR^1$;
\item $\Xi^1_{-1}$ is affinely equivalent to any triangle $\Delta=\{a,b,c\}$ in $\IR^2$;
\item $\Xi^2_0$ is affinely equivalent to an octahedron $\Delta\cup (x-\Delta)$ where $\Delta\subset\IR^2$ is a triangle centered at zero and $x\in\IR^3$ does not belong to the linear hull of the triangle $\Delta$;
\item $\Xi^3_1$ is affinely equivalent to $(x-\Delta)\cup\Delta\cup (x+\Delta)$ where $\Delta=\{a,b,c,d\}$ is a tetrahedron centered in zero and $x\in\IR^4$ does not belong to the linear hull of $\Delta$.
\end{itemize}

In \cite{Ba98} it was proved that for $r\le 3$ any two $r$-centerpole subsets of $\IR^r$ of cardinality $\cp_r^0(\IR^r)$ are affinely equivalent.

\begin{problem} Is each 12-element 4-centerpole set in $\IR^4$ affinely equivalent to the $\binom{3}{1}$-sandwich $\Xi^3_1$?
\end{problem}

\begin{problem} Calculate the (5,0)-centerpole number $\cp^0_5(\IR^5)$ of the Euclidean space $\IR^5$. Theorem~\ref{t10.3} yields the bounds $18\le \cp_5^0(\IR^5)\le \cp_5^0(\IZ^5)\le 25$.
\end{problem}

Observe that the upper bound for $\cp_r^0(\IR^r)$ in Theorem~\ref{t10.3} has exponential growth while the lower bound has polynomial growth.

\begin{problem} What is the growth of the sequence $\big(\cp_r^0(\IR^r)\big)_{r=1}^\infty$? Exponential or polynomial?
\end{problem}

The proof of Theorem~\ref{t9.1} implies that $\cp_r^0(\IR^n)\le\cp_r^\lambda(\IR^n)<2^n$ for any finite nonnegative number $r,\lambda,n$ with $r+\lambda\le n+1$. So, we have some lower and upper bounds for the centerpole numbers $\cp_r^\lambda(\IR^n)$. It would be nice to improve these bounds or find exact values of these centerpole numbers.

\begin{problem} Calculate or evaluate the centerpole numbers $\cp_r^\lambda(\IR^n)$ for various numbers $r,\lambda$ and $n$.
\end{problem}

Concerning the centerpole numbers $\cp_r^\lambda(\IH^2)$ of the hyperbolic plane, we know their values only for $r\le 2$:

\begin{theorem} $\cp_1^\lambda(\IH^2)=1$ and  $\cp_2^\lambda(\IH^2)=3$ for any $\lambda\in\{0,1,2\}$.
\end{theorem}

For $r>2$ we have some information on the centerpole numbers $\cp_r^0(\IH^n,\asdim|\Borel,\mathcal S)$ of the $n$-dimensional hyperbolic space $\IH^n$ endowed with the measure $\asdim|\Borel$ defined on the $\sigma$-algebra $\Borel$ of Borel subsets of $\IH^n$. This centerpole number will be denoted by $\Bcp_r^0(\IH^n)$. Observe that $\Bcp_r^0(\IH^n)$ is equal to the smallest cardinality $|\mathcal S'|$ of a subfamily $\mathcal S'\subset\mathcal S=\{s_c:c\in\IH^n\}$ such that for each Borel $r$-coloring of $\IH^2$ there is an unbounded monochromatic $\mathcal S'$-symmetric subset $A\subset \IH^n$.

\begin{theorem} $3\le \Bcp_r^0(\IH^n)\le\non(\M)$ for any numbers $r,n\ge 2$.
\end{theorem}

Here $\non(\M)$ denotes the smallest cardinality of a non-meager subset of the real line $\IR$. It is known that $\aleph_1\le \non(\M)\le\mathfrak c$ and the position of the cardinal $\non(\M)$ on the interval $[\aleph_1,\mathfrak c]$ depends on Set-Theoretic Axioms, see \cite{Vaughan}.

\section{Dimension symmetry and centerpole numbers of groups}

In this section we study the dimension symmetry numbers $\mds_r(G)$ and centerpole numbers $c_r^\lambda(G)$ of groups $G$, considered as symmetry measure spaces $(G,\asdim,\mathcal S)$ endowed with the family $\mathcal S=\{s_c:c\in G\}$ of central symmetries $s_c:x\mapsto cx^{-1}x$, and the (non-additive) measure $\asdim:\mathcal P(G)\to[-1,\infty]$ assigning to each subset $A\subset G$ its asymptotic dimension $\asdim(A)$, defined as follows.

We say that a subset $A$ of a group $G$ has asymptotic dimension $\asdim(A)\le n$ for some $n\ge-1$ if for each finite subset $F\subset G$ there is a cover $\U$ of $A$ with finite $\mesh(\U)=\bigcup_{U\in\U}U^{-1}U$ such that for some finite subset $B\subset G$ and each point $x\in G\setminus B$ the set $xF$ meets at most $n+1$ sets of the family $\U$.

We write $\asdim(A)=n$ if $\asdim(A)\le n$ but $\asdim(A)\not\le n-1$. If $\asdim(A)\not\le n$ for all $n\in\IN$, then we write $\asdim(A)=\infty$. This definition implies that a subset $A\subset G$ is finite if and only if $\asdim(A)=-1$.

To shorten notations, by $\mds_r(G)$ we shall denote the maximal symmetry number $\ms_r(G,\asdim,\mathcal S)$ of the symmetry measure space $(G,\asdim,\mathcal S)$ and call it the {\em dimension symmetry number} of $G$. Observe that $\mds_r(G)\ge 0$ if and only if $\MS^+_r(G)\ge\aleph_1$. Because of that, Theorem~\ref{t8.5} imply the following characterization:

\begin{theorem} An infinite Abelian group $G$ has dimension symmetry number $\mds_r(G)\ge0$ for some cardinal $r\ge 1$ if and only if one of the following conditions is satisfied:
\begin{enumerate}
\item $r\le r_0(G)$ and the group $G$ is finitely-generated;
\item $r\le r_0(G)+1$ and the group $G$ is countable, infinite-generated, and has finite 2-rank $r_2(G)$;
\item $r\le\max\{\log|G|,r_2(G)\}$ and the group $G$ is uncountable or has infinite 2-rank $r_2(G)$.
\end{enumerate}
\end{theorem}

In the same way we can reformulate Theorems~\ref{t8.15}--\ref{t8.17} of Yu.~Gryshko and A.~Khelif.

\begin{problem} Given an (Abelian) group $G$ calculate or evaluate its dimension symmetry numbers $\mds_r(G)$.
\end{problem}

A similar problem can be asked for the centerpole numbers $\cp_r^\lambda(G)$ of the symmetry measure spaces $(G,\asdim,\mathcal S)$.

\begin{problem} Given a group $G$ calculate or evaluate its centerpole numbers $\cp_r^\lambda(G)$.
\end{problem}

Some information on the centerpole numbers $\cp_r^0(G)$ of abelian groups $G$ can be derived from Theorems~\ref{t10.2} and \ref{t10.3}. Much less is known about the centerpole numbers of non-commutative groups, in particular, of the free group $F_2$ with two generators.

\begin{problem} Has the free group $F_2$ finite centerpole numbers $c_r^0(F_2)$ for all $r\in\IN$?
\end{problem}

For $r=2$ the answer to this problem is affirmative. Analyzing the proof of Theorem~1 of \cite{GK05}, one can see that $c_2^0(F_2)\le 4$.

\begin{problem} Is $c_2^0(F_2)$ equal to 4, 3, or 2?
\end{problem}

\end{document}